\documentclass[a4paper]{article}

\usepackage{MyLaTeX,MyArticle}

\newcommand{\Vmin}{V_{\min}}
\title{Subspace stabilizers and maximal subgroups of exceptional groups of Lie type}
\author{David A.\ Craven}
\date{\today}

\begin{document}
\maketitle

\begin{abstract}
In 1998, Liebeck and Seitz introduced a constant $t(G)$, dependent on the root system of a reductive algebraic group $G$ and proved that if $x$ is a semisimple element of order greater than $t(G)$ in $G$ then there exists an infinite subgroup of $G$ stabilizing the same subspaces of $L(G)$ as $x$. The values for $t(G)$ are $12$, $68$, $124$ and $388$ for $G=G_2,F_4,E_6,E_7$ respectively. In this paper we obtain a similar result for these groups and the minimal module $\Vmin$, obtaining significantly smaller numbers, namely $4$, $18$, $27$ and $75$ respectively (with some small conditions on the element $x$ that are not important for applications). Note that both $t(G)$ and these new bounds are sharp. As a corollary we eliminate several potential maximal subgroups $\PSL_2(q_0)$ of these groups that seem difficult to eliminate through other means, along with other groups.

This paper forms part of the author's programme to vastly reduce the number of putative maximal subgroups of exceptional groups of Lie type.
\end{abstract}
\section{Introduction}

The determination of the maximal subgroups of the exceptional groups of Lie type is one of the main outstanding problems in the theory of the finite simple groups. All remaining unclassified maximal subgroups are almost simple groups, and the bulk of that list are Lie type groups in the same characteristic, and of these the list is dominated by $\PSL_2(q)$.

In \cite{liebeckseitz1998}, Liebeck and Seitz prove a theorem that states that, for each reductive algebraic group $G$ there is an integer $t(G)$, defined in terms of the root system of $G$, such that, if $x$ is a semisimple element of order greater than $t(G)$ in $G$, there exists an infinite subgroup $X$ of $G$ such that $x$ and $X$ stabilize the same subspaces of the adjoint module $L(G)$.

They then apply this to get bounds on the possible $q_0$ such that $\PSL_2(q_0)$ can be a maximal subgroup of a finite group of Lie type $G(q)$ (where $q$ and $q_0$ are both powers of the same prime $p$). If $p$ is odd, we get $q_0<2 t(G)$, and if $p=2$ we get $q_0<t(G)$.

The purpose of this paper is to achieve a similar result for the minimal module $\Vmin$, but with much smaller bounds than $t(G)$, thus reducing significantly the possible $q_0$ such that an almost simple group with socle $\PSL_2(q_0)$ can be a maximal subgroup of an exceptional group of Lie type (unless it is already known).

\begin{thm}\label{thm:mainthm}
Let $G$ be one of the groups $G_2$, $F_4$, $E_6$ and $E_7$ in characteristic $p\geq 0$, and let $\Vmin$ denote the minimal module of dimension $7-\delta_{p,2}$, $26-\delta_{p,3}$, $27$ and $56$. If $x$ is a semisimple element of order $n$ in $G$, then there exists a positive-dimensional subgroup $X$ of $G$ stabilizing the same subspaces of $\Vmin$ as $x$ if one of the following holds:
\begin{enumerate}
\item\label{thmp:g2} $G=G_2$ and $n\geq 5$;
\item\label{thmp:f4} $G=F_4$ and $n\geq 19$;
\item\label{thmp:e6real} $G=E_6$, $x$ is real, and $n\geq 19$;
\item\label{thmp:e6n28} $G=E_6$, and $n\geq 28$ with $6\nmid n$ and $Z(G)\not\leq \gen x$;
\item\label{thmp:e6odd} $G=E_6$ and $n\geq 77$ with $n$ odd;
\item\label{thmp:e7}$G=E_7$ and $n\geq 77$ with $n$ odd.
\end{enumerate}
\end{thm}

The bound $t(G)$ for these groups is $12$, $68$, $124$ and $388$ respectively so these are significant improvements. (Note that these bounds are tight, in the sense that there are elements of order $t(G)$ in $G$ that do not have the required property.) The caveats on the congruences of the orders of $x$ are artefacts of the particular representation of semisimple elements that we choose inside $G$, and could be removed with a different choice, which we do for $F_4$ but not for larger groups; however, this version is sufficient for our needs.

As a consequence of this, we can eliminate various maximal subgroups of exceptional groups of Lie type, particularly Lie type groups in the same characteristic as $G$, as was done in \cite{liebeckseitz1998} using their weaker version of Theorem \ref{thm:mainthm}. If $H=\PSL_2(q_0)$ (or possibly $\SL_2(q_0)$ in the case of $E_7$, since the group has a centre of order $2$ when acting on the minimal module) is a subgroup of $G$ then $H$ has elements of order $(q_0\pm 1)/2$, one of which is guaranteed to be odd. We therefore get the following corollary.

\begin{cor}\label{cor:psl2} Let $H$ be a maximal subgroup of the exceptional simple group of Lie type $G(q)$, not of type $E_8$, with $F^*(H)=\PSL_2(q_0)$. If $H$ is not the fixed points of an $A_1$ subgroup of the algebraic group $G$ (and hence known) then one of the following holds:
\begin{enumerate}
\item $G=F_4(q),E_6(q),{}^2\!E_6(q)$ and $q_0\leq \gcd(2,q_0-1)\cdot 18$;
\item $G=2\cdot E_7(q)$, and $q_0\leq \gcd(2,q_0-1)\cdot 75$;
\end{enumerate}
\end{cor}

This includes, as just one example, the case of $\PSL_2(128)\leq E_7(2^a)$, which is very difficult to remove via other means as there are few semisimple elements of small orders to use to get restrictions on the possible composition factors of the restrictions of the minimal and adjoint modules. In subsequent work the author will address these smaller $\PSL_2$s, but the techniques there will be more akin to those used in \cite{craven2015un2} and less like those employed here. Theorem \ref{thm:mainthm} can be used to eliminate still other cases left open in \cite{liebeckseitz1998}, but we do not list them explicitly here.
%. Recall from \cite{liebeckseitz1998} that the remaining cases are Lie type groups $H(q_0)$ up to untwisted rank half that of $G$, with $q_0\leq 9$, and $\PSL_3(16)$ and $\PSU_3(16)$.

%\begin{cor}
%Let $G=G(q)$ be an exceptional group of Lie type not of type $E_8$, and let $H=H(q_0)$ be a maximal subgroup of $G$, with $q$ and $q_0$ being powers of a prime $p$. Suppose that $H$ is not the fixed points of a positive-dimensional subgroup of the corresponding algebraic group $\mathbb{G}$ under a Frobenius endomorphism of $\mathbb{G}$. Then one of the following holds:
%\begin{enumerate}
%\item $F^*(H)$ is $\PSL_2(q_0)$, a small Ree or Suzuki group, and the possibilities are given above.
%
%\item $F^*(H)=\PSL_3(q_0)$, $q_0=2,3,4$ ($G=F_4$), $q_0=2,3,4,7$ ($G=E_6$), $q_0=2,3,4,5,7,8$ ($G=E_7$)
%\item $F^*(H)=\PSU_3(q_0)$, $q_0=3,4,5$ ($G=F_4$), $q_0=3,4,5,8$ ($G=E_6$), $q_0=3,4,5,7,8$ ($G=E_7$)
%\item $F^*(H)=\PSp_4(q_0)'$, $q_0\leq 5$ ($G=F_4,E_6$), $q_0\leq 9$ ($G=E_7$)
%\item $F^*(H)=G_2(q_0)'$, $q_0=2,3$ ($G=F_4,E_6$), $q_0\leq 8$ ($G=E_7$)
%\item $F^*(H)=\PSL_4(q_0)$, $q_0=2,3$ ($G=E_6$) $q_0=2,3,5$ ($G=E_7$)
%\item $F^*(H)=\PSU_4(q_0)$, $q_0=2,3$ ($G=E_6$) $q_0=2,3,4,5$ ($G=E_7$)
%\item $F^*(H)=\PSp_6(q_0)$, $q_0=2$ ($G=E_6$), $q_0=2,3$ ($G=E_7$)
%\item $F^*(H)=\POmega_7(q_0)$, $q_0=2$ ($G=E_6$), $q_0=2,3$ ($G=E_7$).
%\end{enumerate}
%\end{cor}
%\begin{pf} This is a simple application of the main theorem, once one computes the orders of elements inside these groups. The semisimple elements of groups of Lie type are well known, or one can use a computer.
%\end{pf}

If one wants maximal subgroups of almost simple groups whose socle is an exceptional group of Lie type, rather than just the simple group, then one needs to be concerned about those outer automorphisms that do not stabilize the minimal module, i.e., graph automorphisms for $F_4$ and $p=2$ and for $E_6$. We discuss this in Section \ref{sec:almostsimple}, in particular extending Corollary \ref{cor:psl2} to such groups.

\medskip

The methods in this paper are both theoretical and computational. The reason for computers being involved is that, in order to understand the bounds for the minimal module $\Vmin$, we need to work with subgroups of a maximal torus, and eventually determine the exponents of many finitely presented abelian groups. While each one would be easy for a human to do, the sheer number of possible outcomes means that a computer is needed. Complete Magma code to generate the relevant groups and determine their exponents is given on the author's website, along with the outputs of the algorithms.

Section \ref{sec:theory} introduces the theoretical underpinnings for statements like Theorem \ref{thm:mainthm}, starting from a general setup. Section \ref{sec:practice} then gives details on the computation for each of $G_2$, $F_4$, $E_6$ and $E_7$, with $G_2$ being done by hand and the others being done by computer. The case of $G_2$ is sufficiently detailed that the reader can see how the computations go in the larger cases. We then use our methods to independently verify the results of Liebeck--Seitz for $L(G)$ for $G=G_2,F_4$. Section \ref{sec:almostsimple} discusses extending Corollary \ref{cor:psl2}, and other such corollaries that depend on Theorem \ref{thm:mainthm} to almost simple groups whose socle is an exceptional group of Lie type.

\section{Theory}
\label{sec:theory}

Let $G$ be a reductive algebraic group over a field $k$ of characteristic $p$, and let $M$ be a $kG$-module of dimension $d$. Let $T$ denote a maximal torus in $G$, and assume that $T$ acts diagonally on $M$. Let $x$ be a semisimple element of $G$ lying in $T$ of order $n$, and let $\zeta$ denote a primitive $n$th root of unity: the $\zeta^i$-eigenspaces are subspaces of $M$ spanned by sets of standard basis vectors $e_j$, since $x$ acts diagonally.

Notice that if $y$ stabilizes all of the subspaces of $M$ that $x$ stabilizes then in particular each $e_i$ is an eigenvector for $y$, so that $y$ acts diagonally and $y\in T$. Thus it makes sense to focus our attention on subgroups of $T$.

\begin{prop}
There exists a finite set of integers $\mc X_M$ such that if $n\notin X_M$ then for any semisimple element $x\in T$ of order $n$ there exists an infinite subgroup $Y$ of $T$ such that $x$ and $Y$ stabilize the same subspaces of $M$, and conversely, if $n\in \mc X_M$ then there exists a choice of $x$ of order $n$ such that there is no infinite subgroup $Y$ stabilizing the same subspaces of $M$ as $x$.
\end{prop}
\begin{pf} If $W$ of dimension $r$ is an eigenspace of $x$ then $W$ is spanned by $\{e_{i_1},e_{i_2},\dots,e_{i_r}\}$, where $e_1,\dots,e_d$ denotes the standard basis for $M$. In particular, there are only finitely many possibilities for $W$, namely $2^d$. Thus if we collect together the elements of $T$ according to their sets of eigenspaces, there are only finitely many such collections. Thus there are only finitely many subgroups $S$ of $T$ that are the stabilizers of these sets of eigenspaces.

Let $\mc X_M$ denote the set of orders of elements of those subgroups $S$ of $T$ that are finite. If $x\in T$ has order not in $\mc X_M$ then $x$ must lie in an infinite subgroup $S$ stabilizing the same subspaces of $M$ as $x$, as needed.

The converse is clear from the definition of $\mc X_M$.
\end{pf}

\begin{rem}\label{rem:notcoprime}
Let $x$ be semisimple of order $n$. Just because $x$ is contained in an infinite subgroup $X$ stabilizing the same subspaces of a module $M$, does not mean that for any positive integer $a$ there is an element $y$ such that $o(y)=an$, $y^a=x$ and $x$ and $y$ stabilize the same subspaces, it merely asserts that there are infinitely many such $a$. For example, using a computer it can be shown that for $G=F_4$ and $n=70$, $a=2$, there are conjugacy classes that fail this for $M=L(G)$, even though $70>t(F_4)$ so this is covered by \cite[Theorem 1]{liebeckseitz1998}.

This appears to only be the case where $a$ and $n$ are not coprime.
\end{rem}

What we see from the proposition is that, although there are many subgroups of $T$, only finitely many of them appear as the stabilizers of the set of subspaces that are stabilized by an element of $T$. If $M=L(G)$ then the result of Liebeck--Seitz referenced in the introduction states that $\mc X_M$ is a subset of $\{1,\dots,t(G)\}$, and their proof shows that $\max \mc X_{L(G)}=t(G)$. We will use $\Vmin$ rather than $L(G)$ in an attempt to get better bounds, at the expense of having to use more effort. However, much of the effort of millions of calculations in abelian groups is done via computer.

We introduce a few pieces of notation and some definitions to help our discussion later. We start with an omnibus definition, containing many of the basic ideas we need.

\begin{defn} Let $e_1,\dots,e_d$ be a basis for $M$ such that $T$ acts diagonally on the $e_i$. A \emph{block system} $B$ is a set partition of $\{1,\dots,d\}$, whose constituent sets are called \emph{blocks}, and the \emph{stabilizer} of $B$ is the subgroup of $T$ that acts as a scalar on every subspace $W$ of $V$ such that $W$ is spanned by the $e_i$ for $i$ in a block of $B$. The \emph{dimension} of a block system $B$, denoted $\dim(B)$, is the dimension of its stabilizer as an algebraic group.

If $B$ and $B'$ are block systems then $B'$ is a \emph{coarsening} of $B$ if any block of $B$ is a subset of a block of $B'$, in other words, if the blocks of $B'$ are obtained by amalgamating blocks of $B$.

If $x$ is an element of $T$, then the \emph{block system associated to $x$} is the block system where $i$ and $j$ lie in the same set if and only if $e_i$ and $e_j$ are eigenvectors with the same eigenvalue under the action of $x$.

If $B$ is a block system then the \emph{closure} of $B$ is the block system $B'$ such that the stabilizers of $B$ and $B'$ are the same and if $B''$ is any coarsening of $B$ such that $B$ and $B''$ have the same stabilizer, then $B''$ is a coarsening of $B'$. In other words, the closure of $B$ is the coarsest block system with the same stabilizer as $B$. 
\end{defn}

We see that $x$ is contained in an infinite subgroup $Y$ such that $x$ and $Y$ stabilize the same subspaces if and only if the block system associated to $x$ has positive dimension.

We therefore wish to construct all block systems with finite stabilizers and compute the exponents of such groups. (Since the stabilizers are abelian, we only need the exponents since a finite abelian group contains an element of order $m$ if and only if $m$ divides the exponent of the group.)

In order to construct finite stabilizers, we need some representation for torus elements acting on the minimal module. While this is easy for classical groups, for exceptional groups it is not necessarily so easy to get a representation of a maximal torus acting on $\Vmin$.

Our solution to this is to choose a maximal-rank subgroup of an exceptional algebraic group that is a product of classical groups, and then use the maximal torus from that, which we understand. One downside to this is that the subgroup that we have, for example $A_2A_2A_2$ inside $E_6$, is not really $\SL_3^{\times 3}$ but a central product of $\SL_3$s with a group of order $3$ on top, and hence our diagonal subgroup is a subgroup of index $3$ inside the maximal torus of $E_6$. This means that we only make statements about elements of order prime to $3$ in this case.

In the next section we use specific embeddings of maximal-rank subgroups to yield the theorems stated in the introduction.

\section{Practice}
\label{sec:practice}

The examples in this section all follow a pattern, which we describe now. It might be of benefit to the reader to read this in tandem with the first example for $G_2$ below, which is worked out by hand. Later examples are too large for this and so we resort to computer calculations.

We take an exceptional algebraic group $G$, take a maximal-rank subgroup $H$ of $G$ that is a product of groups of type $A$, and consider the restriction of $\Vmin$ to $H$. We note the order $m$ of $H/H'$ (e.g., it has order $3$ in the case of $A_2A_2A_2\leq E_6$, order $2$ in the case of $A_7\leq E_7$, and so on). We then take a maximal torus $T$ inside $H$ and consider the representation of $T$ on $\Vmin$. If $x$ is an element of $T\leq H$ of order $n$, and $\zeta$ is a primitive $n$th root of unity, then $x$ can be written as a diagonal matrix with entries $\zeta^{a_1},\dots,\zeta^{a_r}$ for some $r$, since $H$ is a product of groups of type $A$, and so there is a natural representation of $x$.

The eigenvalues of $x$ on $\Vmin$ are powers of $\zeta$, and these can be written as $\zeta^{f_i(a_1,\dots,a_r)}$, where the $a_j$ are as above and the $f_i$ are linear functions depending on the representation of $H$ on $\Vmin$.

We then consider the free abelian group on $r$ variables, subject to various relations. The first relations are that the sums of various $a_i$ are zero, coming from the fact that the type-$A$ groups are $\SL_n$ rather than $\GL_n$. Other relations arise from setting various eigenvalues equal to one another. The symmetry group acting on the eigenvalues can be used to reduce the number of possibilities here.

Given a block system, we take all coarsenings of it, take their closures, and consider all their images under the group of symmetries of the eigenvalues, slowly building up a list of all possible stabilizers. We continue doing this until we reach a group with torsion-free rank zero, and take its exponent. We then collate all such exponents, and multiply them by the integer $m$ above, since this is the index $|H:H'|$, and we have only calculated the exponents for $T\cap H'$. The set $\mc X_{\Vmin}$ is a subset of all divisors of these integers. (It would be equal to $\mc X_{\Vmin}$ if $m=1$, but if $m\geq 1$ there is ambiguity in whether certain numbers really lie in $\mc X_{\Vmin}$ or not.)

\subsection{$G_2$}

We start with $G_2$, where we let $H$ be the maximal-rank $A_2$ subgroup. The representation of this on $\Vmin$ has composition factors $L(01)$, $L(10)$ and $L(00)$.

Let $x$ be an element of order $n$ in $A_2$, with eigenvalues $\zeta^a,\zeta^b,\zeta^{-a-b}$, where $\zeta$ is a primitive $n$th root of unity. The eigenvalues of $x$ on $\Vmin$ are
\[ \zeta^a,\zeta^b,\zeta^{-a-b},\zeta^{-a},\zeta^{-b},\zeta^{a+b},1.\]
(We label the basis elements $e_1$ to $e_7$ in the order above.) The block system associated to $x$ is, generically, the singleton partition of $\{1,\dots,7\}$, and of course has dimension $2$. By considering just the exponents of the eigenvalues we get the list
\[ a,b,(-a-b),-a,-b,(a+b),0,\]
and equalities between these yield systems of linear equations. Up to the $\Sym(3)$ automorphism group acting on the torus of $\SL_3$, we may assume that in any refinement $B$ the first element does not lie in a singleton set. We conclude that $a$ equals one of $b$, $-a$, $-b$ or $0$.

Suppose firstly that one of the first six eigenvalues is equal to $1$, say $\zeta^a$. The eigenvalues exponents are therefore $0,b,-b,0,b,-b,0$, but the dimension of the block system
\[ \left\{\bitbig\{1,4,7\},\{2,5\},\{3,6\}\right\}\]
is still $1$, so we need to make more eigenvalues equal. This means that $b$ is equal to either $b$ or $0$, yielding $o(x)=1,2$.

We therefore assume that the $1$-eigenspace of $x$ is $1$-dimensional, i.e., none of $a,b,a+b$ is equal to $0$ in the stabilizer of the block system. We will remove the $1$-eigenspace from our lists from now on to remind us that it has been eliminated.

We still have that $a$ is equal to another eigenvalue exponent, say $b$. In this case the eigenvalue exponents are
\[ a,a,-2a,-a,-a,2a,\]
and the corresponding block system $\{\{1,2\},\{4,5\},\{3\},\{6\},\{7\}\}$ has dimension $1$ again. Setting two of $\pm a$ and $\pm 2a$ equal to one another yields $\alpha a=0$ for some $\alpha=1,2,3,4$, so $o(x)=1,2,3,4$.

We may therefore assume that no two of $a$, $b$ and $-a-b$ are equal to each other, by applying the $\Sym(3)$ automorphism group. This means that $a$ must be equal, via the automorphisms, to either $-a$ or $-b$.

If $a=-b$ then $\zeta^{a+b}=1$, not allowed since we already assume that the $1$-eigenspace is $1$-dimensional. Thus $a=-a$, so $\zeta^a=-1$ (as $\zeta^a\neq 1$). Thus the eigenvalues are
\[ -1,\zeta^b,-\zeta^b,-1,\zeta^{-b},-\zeta^{-b},1.\]
Setting $\zeta^b$ equal to $-1$ gives $o(x)=2$, to $-\zeta^b$ is impossible, to $\zeta^{-b}$ gives $o(x)=2$, and to $-\zeta^{-b}$ gives $\zeta^b=\pm\mathrm{i}$, so that $o(x)=4$.

We therefore have proved the following proposition, giving Theorem \ref{thm:mainthm}(i).

\begin{prop}
If $G=G_2$ and $M=\Vmin=L(10)$ then $\mc X_M=\{1,2,3,4\}$. Consequently, if $H$ is a subgroup of $G_2$ in any characteristic, and it contains a semisimple element of order at least $5$, then $H$ is contained in a positive-dimensional subgroup of $G$ stabilizing the same subspaces of $\Vmin$.
\end{prop}

Comparing this bound with the Liebeck--Seitz bound of $12$, it is clear that we are making significant progress, and it is worth attempting to move on to larger exceptional groups.

\subsection{$F_4$}

For $G=F_4$, as with $G_2$, we choose a maximal-rank subgroup $H$ in which we can easily represent a torus. In this case we choose two different options: $A_2\tilde A_2$ and $A_1^4$. The first will allow us to consider semisimple elements of order prime to $3$ and the second allows us to consider semisimple elements of odd order.

We let $H=A_1^4$ first: the representation of the this is the sum of two trivials and all six possible ways of tensoring two natural modules and two trivial modules for the four $A_1$s, which means that the exponents for the eigenvalues of a general element in the torus of $H$ is
\[ \{\pm a_i\pm a_j:1\leq i<j\leq 4\}\cup\{0\},\]
with $0$ appearing twice, although this is not important for considering coincidences of eigenvalues. The symmetry group we apply here is $\Sym(4)$ acting in the obvious way on the $A_1$s.

We use a computer to analyse this situation, finding all coarsenings of the singleton system of dimension $4$, taking their closures, working up to automorphism, and continuing until we only have block systems of dimension $0$, finding $1264$ distinct block systems, with the following dimensions:
\[ 4^2,\;\; 3^{11},\;\; 2^{113},\;\;1^{538},\;\;0^{600}.\]
The exponents of the six-hundred torsion subgroups (i.e., stabilizers of block systems of dimension $0$) are all even numbers between $2$ and $36$, so if $x$ is an element of odd order at least $19$ then $x$ can only be contained in infinite $\Vmin$-inequivalent subgroups of a maximal torus, as needed.

\medskip

For $H=A_2\tilde A_2$, the representation of $H$ on $\Vmin$ is as the sum of three modules: the tensor product of the two naturals, the tensor product of the two duals, and the trivial for the $A_2$ by the adjoint representation $L(11)$ for the $\tilde A_2$. We can more easily write down the exponents of the eigenvalues in terms of six variables
\[ \{a_i+a_j:1\leq i\leq 3,\;4\leq j\leq 6\}\cup \{a_i+a_j:1\leq i\leq 3,\;4\leq j\leq 6\}\cup \{a_i-a_j:4\leq i\neq j\leq 6\}\cup\{0\},\]
where $a_1+a_2+a_3=0$ and $a_4+a_5+a_6=0$. This time we have $\Sym(3)\times \Sym(3)$ acting by permuting the $a_i$ for $\{1,2,3\}$ and $\{4,5,6\}$.

We again use a computer to analyse this in the same way as before, finding $9278$ distinct block systems up to automorphism, with the following dimensions:
\[ 4,\;\; 3^{17},\;\; 2^{255},\;\;1^{2123},\;\;0^{6882}.\]
The exponents of the nearly seven thousand torsion subgroups are all multiples of $3$ between $3$ and $54$, so if $x$ is an element of order prime to $3$ and at least $19$ then $x$ is contained in an infinite subgroup of the torus stabilizing the same subspaces of $\Vmin$, as needed.

Putting these two results together we get the same result for any element of order at least $19$ that is not a multiple of $6$, and also for any element of order at least $73$. As $21$ is a good number, so is $42$, as $x$ lies in an infinite subgroup stabilizing the same subspaces of $\Vmin$ whenever $x^i$ does. This idea eliminates more numbers, and we are left with only
\[ 24,\;\;30,\;\;36,\;\;48,\;\;72.\]
We check by a computer program that for every one of the $584$ conjugacy classes of element of order $24$, there exists one of order $120$ that powers to it and stabilizes the same subspaces. Therefore elements of order $24$, $48$ and $72$ (the latter two being powers of $24$) are also good, leaving just $30$ and $36$. We do the same thing for the $2333$ classes of elements of order $36$, finding elements of order $180$ that power to them and stabilize the same subspaces of $\Vmin$.

We might want to do the same thing for $30$, but we cannot guarantee that $150$ will work (as $5$ is not prime to $30$, see Remark \ref{rem:notcoprime}), and indeed there are five conjugacy classes that do not have elements of order $150$ that power to them and fix the same subspaces of $\Vmin$. For those remaining classes we have to switch to another element order, and the obvious candidate, bearing in mind Remark \ref{rem:notcoprime}, is $210=7\cdot 30$. If one does this, then those remaining five classes, and the remaining 1229 other classes, also work, and so we get Theorem \ref{thm:mainthm}(ii).

\begin{prop}
If $G=F_4$ and $M=\Vmin=L(\lambda_1)$ then $\mc X_M=\{1,\dots,18\}$. Consequently, if $H$ is a subgroup of $F_4$ in any characteristic, and it contains a semisimple element of order at least $19$, then $H$ is contained in a positive-dimensional subgroup of $G$ stabilizing the same subspaces of $\Vmin$.
\end{prop}

%Doing the same thing for the $1234$ classes of elements of order $30$ and the $2333$ classes of elements of order $36$ (finding elements of orders $150$ and $180$ respectively) we eliminate the final two outstanding element orders, leaving us with the result in Theorem \ref{thm:mainthm} that any element of order at least $19$ lies in an infinite subgroup stabilizing the same subspace of $\Vmin$.

\subsection{$E_6$}

Let $G=E_6$. If $x$ is a real semisimple element then $x$ lies inside $F_4$, and the eigenspaces of $x$ on $\Vmin$ are the same as for $F_4$, except the $0$-eigenspace has dimension one greater than for $F_4$. This does not affect the calculations in the previous subsection, and so we get that if $x$ in $E_6$ is real semisimple and has order at least $19$, then $x$ is contained in an infinite subgroup stabilizing the same subspaces of $\Vmin$, Theorem \ref{thm:mainthm}(iii).

The rest of this subsection considers the case where $x$ is non-real. As with $F_4$, we consider two maximal-rank subgroups $H$ for elements of order prime to $2$ and $3$, settling on $A_5A_1$ and $A_2A_2A_2$ respectively.

We start with $A_5A_1$. The action of this subgroup on $\Vmin$ is as two composition factors, $(L(\lambda_4),L(0))$ and $(L(\lambda_1),L(1))$. As with $F_4$, we label the eigenvalue exponents for the $A_5$ factor by $a_1,\dots,a_6$ such that $\sum_{i=1}^6 a_i=0$, and the other $A_1$ as $\pm a_7$. The eigenvalue exponents of $A_5A_1$ on $\Vmin$ therefore become
\[ \{-(a_i+a_j):1\leq i<j\leq 6\}\cup\{a_i\pm a_7:1\leq i\leq 6\}.\]
Now we have $\Sym(6)$ acting by permuting the $a_i$ for $\{1,\dots,6\}$.

We again use a computer to analyse this, finding $33365$ distinct block systems up to automorphism, with the following dimensions:
\[ 6,\;\;5^7,\;\;4^{68},\;\;3^{630},\;\;2^{4154},\;\;1^{12488},\;\;0^{16017}.\]
When considering the exponents of the finite such groups, we need to consider the centre of $G$, which of course will appear in any one of these subgroups. Firstly, the exponents of them are all multiples of $6$ from $6$ to $156$, so this means that if $x$ has order an odd number at least $77$ or prime to $6$ and at least $29$ then $x$ is contained inside an infinite subgroup stabilizing the same subspaces of $\Vmin$, yielding Theorem \ref{thm:mainthm}(v). However, we are also interested in elements of odd order that do not power to a non-identity central element: checking this is also easy, and we find that the possible orders of such elements are all odd numbers between $1$ and $27$, so we also get Theorem \ref{thm:mainthm}(iv) for $n$ odd.

Now we turn to $H=A_2A_2A_2$. The action of this subgroup on $\Vmin$ is as three composition factors, $(10,01,00)$, $(00,10,01)$ and $(01,00,10)$. We label the eigenvalues exponents by $a_i$ for $i=1,\dots,9$, such that the sums $a_1+a_2+a_3$, $a_4+a_5+a_6$ and $a_7+a_8+a_9$ are all zero. The eigenvalue exponents of $H$ on $\Vmin$ are therefore
\[ \{a_i-a_j:1\leq i\leq 3,\;4\leq j\leq 6\}\cup\{a_i-a_j:4\leq i\leq 6,\;7\leq j\leq 9\}\cup \{a_i-a_j:7\leq i\leq 9,\;1\leq j\leq 3\}.\]
In this case, $\Sym(3)\wr \Sym(3)\leq \Sym(9)$ acts by preserving the set partition $\{\{1,2,3\},\{4,5,6\},\{7,8,9\}\}$.

The computer now finds $26498$ distinct block systems, with the following dimensions:
\[ 6,\;\;5^4,\;\;4^{51},\;\;3^{565},\;\;2^{4002},\;\;1^{12162},\;\;0^{9713}.\]
The exponents of the torsion groups are all multiples of $3$ between $3$ and $81$, so if $x$ is a semisimple element of order at least $28$ and order prime to $3$, then $x$ lies in an infinite subgroup stabilizing the same subspaces of $\Vmin$, completing the proof of Theorem \ref{thm:mainthm}(iv), and yielding the following proposition.

%If $o(x)$ is prime to $3$ then certainly $Z(E_6)\not\leq \langle x\rangle$, so we combine the two statements for the two tori to get the following proposition.

\begin{prop}
Let $G=E_6$ and $M=\Vmin=L(\lambda_1)$. If $x$ is a semisimple element of order at least $28$ and not divisible by $6$, and such that $Z(G)\not\leq \langle x\rangle$, then there exists an infinite subgroup of $G$ stabilizing the same subspaces of $\Vmin$ as $x$. Consequently, if $H$ is a subgroup of $E_6$ in any characteristic not containing $Z(G)$, and it contains a semisimple element of order at least $28$ and not divisible by $6$, then $H$ is contained in a positive-dimensional subgroup of $G$ stabilizing the same subspaces of $\Vmin$.
\end{prop}

It is significantly more difficult to compute the conjugacy classes of semisimple elements of large order in $E_6$ than in $F_4$, and so we do not improve this result in this way. Indeed, non-real elements of $E_6$ of large order are not of so much use, and this is good enough for most purposes.

\subsection{$E_7$}

Let $G=E_7$. Since $\Vmin$ has dimension $56$ and the torus has rank $7$, one expects the number of $\Vmin$-inequivalent subgroups to be much higher, and for the programs to take much longer to run, which is true. It also might mean that there are too many block systems to store them all, and so we have to alter our algorithm for computing these slightly.

We let $H$ be the maximal-rank $A_7$ subgroup. The representation of this on $\Vmin$ has composition factors $L(\lambda_2)$ and $L(\lambda_6)$, i.e., the exterior square of the natural and its dual. The exact form of $H$ is $4\cdot \PSL_8\cdot 2$, and so we will only attempt to obtain information about odd-order elements: this also bypasses questions about whether the element $x$ powers to the central involution.

Let $x$ be an element of order $n$ in $H$, with eigenvalues $\zeta^{a_i}$ for $i=1,\dots,8$ with $\sum a_i=0$, where $\zeta$ is a primitive $n$th root of unity. The eigenvalues of $x$ on $\Vmin$ are therefore $\{\zeta^{\pm(a_i+a_j)}\,:\,1\leq i<j\leq 8\}$. We again use a computer to analyse this situation, although we stop when we reach the block systems of dimension $2$, yielding
\[ 7,\;\; 6^5,\;\;5^{47},\;\;4^{626},\;\; 3^{9781},\;\; 2^{116170}.\]
We obviously don't want to try to store the likely million block systems of rank $1$, so for each block system of dimension $2$ we find all coarsenings, and repeat the process until we reach block systems of dimension $0$. This of course introduces computational repetition but reduces the space requirement.

Doing this produces the set of exponents all multiples of $8$ up to $264$, and all $a\equiv 4\bmod 8$ up to $300$. In particular, the odd divisors of these numbers are all (odd) integers up to $75$, and this means that we get the following proposition and Theorem \ref{thm:mainthm}(vi).

\begin{prop}
Let $G=E_7$ and $M=\Vmin=L(\lambda_1)$. The odd elements of $\mc X_M$ are $\{1,3,\dots,75\}$. Consequently, if $H$ is a subgroup of $E_7$ in any characteristic, and it contains a semisimple element of odd order at least $77$, then $H$ is contained in a positive-dimensional subgroup of $G$ stabilizing the same subspaces of $\Vmin$.
\end{prop}

\subsection{$\mc X_{L(G)}$}

As mentioned in the introduction, the bound $t(G)$ is tight in the sense that there is a semisimple element of order $t(G)$ that does not lie in an infinite subgroup stabilizing the same subspaces of $L(G)$. This follows from the proof of \cite[Proposition 2]{liebeckseitz1998}, but it offers a independent check as to correctness of the methods and programs used in this paper.

Indeed, performing the same calculations with $L(G)$ instead of $\Vmin$ yields the set
\[ \{1,2,3,4,5,6,7,8,9,10,12\}\]
for $G_2$ and
\[ \{1,\dots,58\}\cup\{60,62,64,66,68\}\]
for $F_4$.

\section{From simple to almost simple}
\label{sec:almostsimple}

The results in Theorem \ref{thm:mainthm} are for the simply connected simple algebraic group, but for the general maximal subgroup problem we need to consider maximal subgroups of almost simple groups whose socle is an exceptional group of Lie type. We will extend Corollary \ref{cor:psl2} to almost simple groups here.

Since the maximal subgroups are known for ${}^2\!B_2$, ${}^2\!G_2$, ${}^4\!D_2$, $G_2$ and ${}^2\!F_4$, we are interested in $F_4$, $E_6$, ${}^2\!E_6$ and $E_7$ here, as we have said nothing about $E_8$ in this article.

Let $A$ be a subgroup of the outer automorphism group $\Out(G)$ of a simple group $G$. The maximal subgroups of an almost simple group $G.A$ consist of those containing $G$, those whose intersection with $G$ is a maximal subgroup of $G$, and \emph{novelty} maximals $H$, which are stable under a group of automorphisms but no subgroup between $H$ and $G$ is stable under that group.

If $\Vmin$ is $A$-stable, which is true unless $A$ involves the graph automorphism, then if $H\leq G$ is contained in a positive-dimensional subgroup of $G$ stabilizing the same subspaces of $\Vmin$ then the same is true of any subgroup containing $H$, even in the almost simple group $H.A$. Thus we focus on the case where $A$ contains a graph automorphism of $G$.

Note that, since $L(G)$ is always $A$-stable, if $H$ contains a semisimple element of order not in $\mc X_{L(G)}$ then $H$ is contained in a positive-dimensional subgroup stabilizing the same subspaces of $L(G)$, even in the almost simple group.

If $G=F_4$ and $p=2$ then there is a graph automorphism. In this case, we need to consider  $\SL_2(2^a)$ for $a\geq 5$ for Corollary \ref{cor:psl2}. Note that $\SL_2(2^a)$ contains an element of order $2^a+1$, which does not lie in $\mc X_{L(F_4)}$ for any $a\geq 6$, so we only need to consider $a=5$. In this case the automorphism group has order $5$, so if $\SL_2(32)$ is stabilized by the graph automorphism then it is centralized by it and so cannot yield a maximal subgroup of $G.A$.

We also must consider $G=E_6$ and all primes. Here the graph automorphism swaps $\Vmin$ and its dual, but if $H=\PSL_2(p^a)$ and $x\in H$ has order $n=(p^a\pm 1)/\gcd(p-1,2)$ then $x$ is real, so $x$ lies inside $F_4$ and centralizes at lest a $3$-space on $\Vmin$, including an $F_4$-point of $\Vmin$ (i.e., a point whose line stabilizer is $F_4$). If $X$ is a positive-dimensional subgroup stabilizing the same subspaces of $\Vmin$ as $x$ then $X$ stabilizes the $F_4$-point (or at least the line containing it) so $X\leq F_4$ since $F_4$ is maximal in $G$. But this means that $X$ acts on $\Vmin$ and $\Vmin^*$ in the same way, since $F_4$ does, and so $x$ and $X$ stabilize the same subspaces of $\Vmin\oplus \Vmin^*$, which \emph{is} $A$-stable. Again, we see no novelty maximal subgroups of $G.A$.

%\medskip
%
%We end with a short lemma that could be of interest, concerning the adjoint module for $E_6$.
%
%\begin{lem} Let $x$ be a real semisimple element of odd order in the simply connected form $G$ of $E_6$. If $x$ has order not in $\mc X_{L(G)}$ then $x$ is contained in an infinite subgroup of $G$ stabilizing the same subspaces of $L(G)$.
%\end{lem}

%$A_4$ inside $E_7$, acts like $0$, $\pm(a_i+a_j)$, $\pm a_i$, so $31$ eigenvalues. Get
%\[ \mc X_M=\{1,\dots,28\}\cup \{30\}.\]

\bigskip

\noindent\textbf{Acknowledgements:} I would like to thank the Royal Society for a University Research Fellowship.

\bibliography{references}

\end{document}